\documentclass[11pt,twoside]{article}
\usepackage{atmp}
\copyrightnotice{2002}{6}{197}{205}

\newcommand{\Complex}{\mathbb{C}}
\newcommand{\Integer}{\mathbb{Z}}
\newcommand{\Natural}{\mathbb{N}}

\newcommand{\proj}{\mc{P}_{\!\!\omega}}
  
\newcommand{\mf}{\mathfrak} 
\newcommand{\mc}{\mathcal} 
\def\asg{{\hat{\mf{g}}}}

\newtheorem{theorem}{Theorem}
\newtheorem{lemma}{Lemma}
\newtheorem{corollary}{Corollary}

\begin{document}

\url{math.QA/0203133}
  
\title{An algorithm for twisted\\fusion rules}

\author{Thomas Quella}
\address{Max-Planck Institut f\"ur Gravitationsphysik\\
  Albert-Einstein-Institut, Am M\"uhlenberg 1\\ D-14476 Golm, Germany}
\addressemail{quella@aei-potsdam.mpg.de}
\author{Ingo Runkel and Christoph Schweigert}
\address{LPTHE -- Univ.\ Paris VI,
  4 place Jussieu\\F -- 75252 Paris, France.}
\addressemail{ingo@lpthe.jussieu.fr\qquad schweige@lpthe.jussieu.fr}

\markboth{\it AN ALGORITHM FOR TWISTED FUSION RULES\ldots}{\it T.\ Quella, I.\ Runkel and Ch.\ Schweigert}
  
\begin{abstract}
We present an algorithm for an efficient calculation of the
fusion rules of twisted representations of untwisted affine Lie algebras.
These fusion rules appear 
in WZW orbifold theories and as annulus coefficients in boundary WZW 
theories; they provide NIM-reps of the WZW fusion rules.
\end{abstract}

\cutpage

It is a well known fact that affine Lie algebras have twisted
integrable highest weight representations, and also their fusion rules 
can be determined \cite{Gaberdiel:1997kf,Birke:1999ik}.
The study of conformal field theories provides two interpretations 
for these algebraic objects: They appear as fusion rules in WZW orbifolds,
and on surfaces with boundaries twisted representations label symmetry 
breaking boundary conditions; their fusion rules describe annulus 
coefficients \cite{Birke:1999ik}, see also
\cite{Fuchs:2000vg,Ishikawa:2001zu,Gaberdiel:2002qa}. In this
short note we present an algorithm to compute these fusion rules 
efficiently.

More precisely, we work in the following setting.
Let $\asg^{(1)}$ be an untwisted affine Lie algebra and $\omega$ an
automorphism of order $N$ of its horizontal subalgebra $\mf{g}$.
In the WZW theory based on $\asg^{(1)}$ at level $k$ with 
modular invariant given by charge conjugation, we consider
boundary conditions for which left movers and right movers are related by 
the automorphism $\omega$ at the boundary. By T-duality, these boundary 
conditions correspond to symmetry preserving boundary conditions in a theory 
with modular invariant of automorphism type $\omega$.
This kind of boundary conditions was analysed in \cite{Birke:1999ik}
and more recently again in \cite{Q,Petkova:2002yj,Gaberdiel:2002qa,AFQS}.

The spectrum of open strings living between two boundary
conditions $\alpha,\beta$ is encoded in the boundary partition function
\begin{eqnarray*}
  Z_{\beta\alpha}(q)
  =\sum_{i}N_{i\alpha}^\beta\:\chi_i(q)
\end{eqnarray*}
where the sum over $i$ runs over integrable highest weight representations
of $\asg^{(1)}$ at level $k$ and $\chi_i(q)$ are the corresponding 
characters. The set of boundary conditions is given by twisted representations 
of $\asg^{(1)}$ at level~$k$ and the annulus coefficients $N_{i\alpha}^\beta$ 
are the corresponding twisted fusion rules~\cite{Birke:1999ik}. 
They form a
representation of the fusion rules of the WZW theory at level $k$ by
matrices with non-negative integer entries, a so-called NIM-rep.

\medskip

In order to describe the set of twisted representations we need to introduce
some notation. We denote the weight lattice of the horizontal
subalgebra $\mf{g}$ by~$L$.\footnote{Notice that $L$, in contrast to frequent
  use in the literature,
refers to the weight lattice and {\em not} to the root lattice. This 
convention will be more economic later.} 
A basis of this lattice are the fundamental weights $\Lambda_{(i)}$.
The Killing form endows $L$ with a bilinear form $(\cdot,\cdot)$, and
on $L$ we have  the action of the Weyl group $W$ which is generated by 
reflections
$s_i(\lambda)=\lambda-2\bigl(\lambda,\alpha_{(i)}\bigr)\alpha_{(i)}/\bigl(\alpha_{(i)},\alpha_{(i)}\bigr)$
at the hyperplanes perpendicular to the simple roots $\alpha_{(i)}$. 
The lattice $L^\vee$ dual to $L$ is the coroot lattice of $\mf{g}$;
a basis are the simple coroots $\alpha_{(i)}^\vee$. The lattices $L$ and 
$L^\vee$ inherit an action of the automorphism $\omega$, which can be 
decomposed into an outer automorphism $\omega_0$ and an inner one $\omega_i$,
$\omega=\omega_i\circ\omega_0$. While the inner automorphism $\omega_i$
can be chosen to be the adjoint action of an element of a Cartan subalgebra
and therefore 
induces a trivial action on $L$ and $L^\vee$, the outer part $\omega_0$
can be chosen to be a diagram automorphism of the Dynkin diagram 
of $\mf{g}$. It acts on the lattices $L$ and $L^\vee$ by the
permutations $\omega_0(\Lambda_{(i)})=\Lambda_{(\omega_0 i)}$ and
$\omega_0\bigl(\alpha_{(i)}^\vee\bigr)=\alpha_{(\omega_0 i)}^\vee$ of 
fundamental weights or simple coroots, respectively. Without loss of
generality we can therefore assume $\omega$ to be a diagram automorphism. 
  The length of
the orbit $\bigl\{\Lambda_{(i)},\omega(\Lambda_{(i)}),\omega^2(\Lambda_{(i)}),
\ldots\bigr\}$ will be denoted by $n_i$. 
We also define the lattice of symmetric weights 
$L_\omega=\bigl\{\mu\in L \,\bigl|\,\omega(\mu)=\mu\bigr\}$ which inherits the scalar 
product from $L$. 

An important ingredient in our algorithm is the subgroup
\cite{Fuchs:1996zr,Fuchs:1996vp} 
\begin{eqnarray*}
  W_\omega=\bigl\{\,w\in W \,\big|\, w\circ\omega=\omega\circ w \,\bigr\}
\end{eqnarray*}
of the Weyl group that commutes with the action of $\omega$. It is a
Coxeter group 
with the following generators $\tilde s_i$:
for orbits of length 1, take $\tilde s_i=s_i$. If $i\neq\omega i$, take the 
product $\tilde s_i=s_i s_{\omega i} \ldots s_{\omega^{n_i-1}i}$. This 
prescription needs to be modified, if the element $A_{i,\omega i}$ of the
Cartan matrix is non-vanishing, which in our situation only happens for the
outer automorphism 
of $A_{2n}$ and the orbit consisting of the two nodes in the middle of the 
Dynkin diagram. In this case, take 
$\tilde s_i = s_i s_{\omega i} s_i= s_{\omega i} s_i s_{\omega i}$.

We also need the orthogonal projection of weight space onto its symmetric
subspace:
  $\proj$ defined
  by $\proj=\frac{1}{N}\bigl(1{+}\omega{+}\cdots{+}\omega^{N-1}\bigr)$, 
  $N$ being the
  order of $\omega$.
For the implementation on a computer, one uses directly the action of 
$\tilde s_i$ on symmetric weights:
\begin{eqnarray}
  \label{apeq:reflection}
  \tilde s_i (\lambda) = \lambda - \frac{2(\lambda, \proj\alpha_{(i)})}
{(\proj\alpha_{(i)},\proj\alpha_{(i)})} \proj\alpha_{(i)}\:.
\end{eqnarray}
While the symmetric weight lattice $L_\omega$ is not invariant under the full 
Weyl group,
it admits an action of $W_\omega$. 

We may also define a symmetric coroot lattice
$(L^\vee)_\omega=\bigl\{\beta\in L^\vee\bigl|\omega(\beta)=\beta\bigr\}$.
Note that $L_\omega$ and $(L^\vee)_\omega$ are not dual to each other.
Instead one finds that the lattice $\bigl((L^\vee)_\omega\bigr)^\vee$ dual to
$(L^\vee)_\omega$ involves fractional symmetric weights.
$\proj$ restricts to a surjective map
from $L$ to $\bigl((L^\vee)_\omega\bigr)^\vee$.

We summarise the 
expressions for the different lattices by comparing to the situation for inner 
automorphisms where just two lattices appear:
\begin{itemize}
\item Weight lattice:
  $L=\bigl\{\sum_i\lambda_i\Lambda_{(i)} \,\bigl|\,\lambda_i\in\Integer\bigr\}$\:.
\item Coroot lattice:
  $L^\vee=\bigl\{\sum_i\beta_i\alpha_{(i)}^\vee \,\bigl|\,\beta_i\in\Integer\bigr\}\subset L$\:.
\end{itemize}
  In addition there are four lattices related to the automorphism $\omega$.
\begin{itemize}
\item Symmetric weight lattice:
  $L_\omega=\bigl\{\sum_i\lambda_i\Lambda_{(i)} \,\bigl|\,
   \lambda_i\in\Integer,\:\lambda_{\omega i}=\lambda_i\bigr\}\subset L$\:.
\item Symmetric coroot lattice:\\
  $(L^\vee)_\omega=\bigl\{\sum_i\beta_i\alpha_{(i)}^\vee \,\bigl|\,
   \beta_i\in\Integer,\:\beta_{\omega i}=\beta_i\bigr\}\subset L^\vee$\:.
\item Fractional symmetric weight lattice:\\
  $\bigl((L^\vee)_\omega\bigr)^\vee=\bigl\{\sum_i\lambda_i\Lambda_{(i)}\,\bigl|\,
   n_i \lambda_i\in\Integer,\:\lambda_{\omega i}=\lambda_i\bigr\}\supset L_\omega$\:.
\item Fractional symmetric coroot lattice:\\
  $(L_\omega)^\vee=\bigl\{\sum_i\beta_i\alpha_{(i)}^\vee \,\bigl|\,
   n_i \beta_i\in \Integer,\:\beta_{\omega i}=
   \beta_i\bigr\}\supset(L^\vee)_\omega$\:.
\end{itemize}
  Recall that the $n_i$ are the orbit lengths of fundamental weights.

  The integrable highest weight modules of $\asg^{(1)}$
  at level $k$ are in one-to-one correspondence with elements in
  $P_k^+=L/(W\ltimes kL^\vee)$.
  The expression $W\ltimes kL^\vee$ is just the 
  decomposition of the affine Weyl
  group into a semi-direct product of the finite Weyl group and the 
  translations with
  respect to the scaled coroot lattice. Alternatively, the affine Weyl group
  is generated by finite Weyl reflections and one additional element,
  a shifted Weyl reflection. The latter is a combination of an elementary
  reflection at the highest root $\theta$ of $\mf{g}$ and a
  translation. This amounts to an orthogonal reflection with respect to
  the hyperplane $(\theta,\cdot)=k$.
  An analogous construction can be performed with respect to the lattices
  $L_\omega$ and $\bigl((L^\vee)_\omega\bigr)^\vee$. This defines the sets
  $S_k^+=L_\omega/\bigl(W_\omega\ltimes k(L^\vee)_\omega\bigr)$ and
  $B_k^+=\bigl((L^\vee)_\omega\bigr)^\vee/\bigl(W_\omega\ltimes k(L_\omega)^\vee\bigr)$.
  While $W_\omega\ltimes k(L^\vee)_\omega$ is generated by $W_\omega$
  and the shifted Weyl reflection at $(\theta,\cdot)=k$, the corresponding
  shifted Weyl reflection for $W_\omega\ltimes k(L_\omega)^\vee$ is at
  the hyperplane $(\theta_\omega,\cdot)=k$. The vector
  $\theta_\omega$ in weight space is defined in Table~\ref{aptb:THR}.
\begin{table}[tb]
  \centerline{\begin{tabular}{c|cccccc}
    $\mf{g}$ & $A_{2n}$ & $A_3$ & $A_{2n+1}$ & $D_4 \text{ (triality)}$ & $D_n$ & $E_6$ \\\hline
    $\theta_\omega$ & $2\bigl(\Lambda_{(1)}+\Lambda_{(2n)}\bigr)$ & $2\Lambda_{(2)}$ & $\Lambda_{(2)}+\Lambda_{(2n)}$ & $\Lambda_{(1)}+\Lambda_{(3)}+\Lambda_{(4)}$ & $2\Lambda_{(1)}$ & $\Lambda_{(1)}+\Lambda_{(5)}$
  \end{tabular}}
  \caption{\label{aptb:THR}The vector $\theta_\omega$ in the labeling conventions of~\cite[p.\ 53]{Kac:1990}.}
\end{table}
For each of the three subsets there is a natural choice of a fundamental
domain.
\begin{itemize}
\item Integrable highest weights\\
  $P_k^+=\bigl\{\lambda=\sum_i \lambda_i\Lambda_{(i)} \,\bigl|\,\lambda_i\in\Natural_0\text{ and }(\theta,\lambda)\leq k\bigr\}$\:.
\item Symmetric integrable highest weights\\
  $S_k^+=\bigl\{\lambda=\sum_i\lambda_i\Lambda_{(i)} \,\bigl|\,\lambda_i\in\Natural_0, (\theta,\lambda)\leq k\text{ and }\lambda_{\omega i}=\lambda_i\bigr\}$\:.
\item 
Boundary labels correspond to twisted highest weight
representations \cite{Birke:1999ik} or, equivalently, to irreducible integrable
highest weight representations of the corresponding twisted Lie algebra. They
are labelled by\\
  $B_k^+=\bigl\{\beta=\sum_i\beta_i\Lambda_{(i)} \,\bigl|\:
  n_i \beta_i\in \Natural_0, (\theta_\omega,\beta)\leq k
  \text{ and }\beta_i=\beta_{\omega i}\bigr\}$\:.
\end{itemize}
There is a distinguished vector $\rho_\omega=\sum_in_i^{-1}\Lambda_{(i)}$ in
the lattice $\bigl((L^\vee)_\omega\bigr)^\vee$ which is a fractional analogue
of the Weyl vector $\rho=\sum_i\Lambda_{(i)}$. 
  We denote by $P_k^{++}$, $S_k^{++}$ and $B_k^{++}$ the subsets obtained from
  $P_k^+$, $S_k^+$ or $B_k^+$ after dropping elements which belong to the
  boundary of the respective Weyl chamber, i.e.\ are left invariant by at least
  one nontrivial element of $W \ltimes k L^\vee$,
  $W_\omega\ltimes k(L^\vee)_\omega$ or
  $W_\omega\ltimes k(L_\omega)^\vee$, respectively.
  It is not difficult to see that there exist identifications
  of the form $P_k^++\rho=P_{k+g^\vee}^{++}$, $S_k^++\rho=S_{k+g^\vee}^{++}$
  and $B_k^++\rho_\omega=B_{k+g^\vee}^{++}$
  where $g^\vee$ is the dual Coxeter number of $\mf{g}$. These are a simple
  consequence of the fact that
  $(\theta,\rho)=(\theta_\omega,\rho_\omega)=g^\vee-1$.

We are now prepared to state our result for the determination of twisted
fusion rules. It is a generalisation of the Racah-Speiser algorithm for
tensor product multiplicities (see e.g.\ \cite{FuchsSchweigert}) and the 
Kac-Walton formula \cite{Kac:1990,Walton:1990sc} (see also \cite{Furlan:1990ce,Fuchs:1990hm})
for ordinary fusion rules.
\begin{theorem}
  \label{apth:main}
  The decomposition of the fusion product
\begin{eqnarray*}
  i\star \alpha=\sum_{\beta\in B_k^+}N_{i\alpha}^\beta\:\beta
\end{eqnarray*}
of an untwisted representation $i\in P_k^+$ of $\asg$ and a twisted
representation $\alpha\in B_k^+$ into twisted representations can be 
obtained by the following algorithm:
\begin{enumerate}
\item Compute the weight system $M_i$, including multiplicities, of the finite 
      dimensional irreducible highest weight representation $i$ of the 
      finite dimensional Lie algebra $\mf{g}$.
\item Use $\proj$ to project the set $M_i$ to the lattice of
  fractional symmetric weights.
\item Add the weight $\alpha$ and the twisted Weyl vector $\rho_\omega$ to
  the resulting weights.
\item Use the reflections~\eqref{apeq:reflection} in $W_\omega$ and the shifted reflection
  at the plane $(\theta_\omega,\cdot)=k+g^\vee$ to map the set
  $\proj M_i+\alpha+\rho_\omega$ to the fundamental domain
  $B_{k+g^\vee}^+$. 
\item Discard weights on the boundary $B_{k+g^\vee}^+\backslash 
    B_{k+g^\vee}^{++}$, i.e. those with at least one vanishing
  entry or scalar product with $\theta_\omega$ equal to $k+g^\vee$.
  Supply each remaining contribution, counting multiplicities, with a sign depending on 
  whether the number of reflections has been even or odd.
\item Subtract the twisted Weyl vector $\rho_\omega$. Adding all contributions
  including the relevant multiplicities and signs gives the fusion
  product.
\end{enumerate}
\end{theorem}

We will split the proof into several steps. First, we summarise some earlier
results which will be important in the sequel. It was shown
in~\cite[(2.57)]{Birke:1999ik} that the twisted fusion coefficients for three
weights $i\in P_k^+$ and $\alpha,\beta\in B_k^+$ are given by the formula
\begin{eqnarray}
  \label{apeq:Annulus}
   N_{i\alpha}^\beta
  =\sum_{\mu\in S_k^+}\frac{\bar{S}_{\beta\mu}^\omega S_{i\mu}S_{\alpha\mu}^\omega}{S_{0\mu}}\:.
\end{eqnarray}
where the matrix $S_{\alpha\mu}^\omega$ is given by \cite[(4.6)]{Birke:1999ik}
\begin{eqnarray}
  \label{apeq:TwSMatrix}
  S_{\alpha\mu}^\omega
  =(\text{phase})\,\bigl|\,L_\omega\,/\,(k{+}g^\vee)(L^\vee)_\omega\bigr|^{-1/2}
    \sum_{w\in W_\omega}\epsilon_\omega(w)e^{-\frac{2\pi i}{k+g^\vee}\bigl(w(\alpha+\rho_\omega),\,\mu+\rho\bigr)}
\end{eqnarray}
(see also~\cite[Theorem 13.9]{Kac:1990}).
Note that it carries two different labels
$\alpha\in B_k^+$ and $\mu\in S_k^+$. The symbol $\epsilon_\omega$ denotes
the sign function of $W_\omega$. As the generators of $W_\omega$
may be products of several generators of $W$, in general
the sign function $\epsilon_\omega$ of $W_\omega$ does not coincide with
the restriction of the sign function $\epsilon$ of $W$ to the subgroup
$W_\omega$.
Using Weyl's character formula, the quotient of S~matrices
$S_{i\mu}/S_{0\mu}$ which appears in~\eqref{apeq:Annulus} may be expressed
as
\begin{eqnarray}
  \label{apeq:Character}
  \frac{S_{i\mu}}{S_{0\mu}}
  =\chi_i\Bigl(-\frac{2\pi i}{k+g^\vee}(\mu+\rho)\Bigr)
  =\sum_{j\in M_i}e^{-\frac{2\pi i}{k+g^\vee}(j,\,\mu+\rho)}
\end{eqnarray}
  where $M_i$ denotes the weight system of the finite dimensional
  highest weight module $i$ of $\mf{g}$
  including the multiplicities. If one inserts the
  expressions~\eqref{apeq:TwSMatrix} and~\eqref{apeq:Character} into the
  definition~\eqref{apeq:Annulus} we may write
\begin{eqnarray}
  \label{apeq:AnnulusAux}
  N_{i\alpha}^\beta
  =\sum_{\mu\in S_k^+}f(\mu+\rho)
  =\sum_{\nu\in S_{k+g^\vee}^{++}}f(\nu)
\end{eqnarray}
  where we used the rule $S_k^++\rho=S_k^{++}$ and defined the function
\begin{multline}
  \label{apeq:fAux}
  f(\nu)
  =\bigl|L_\omega/(k+g^\vee)(L^\vee)_\omega\bigr|^{-1}\\\times
   \sum_{j\in M_i}\sum_{w_1,w_2\in W_\omega}\epsilon_\omega(w_1)\epsilon_\omega(w_2)
   e^{-\frac{2\pi i}{k+g^\vee}
   \bigl(\proj j+w_1(\alpha+\rho_\omega)-w_2(\beta+\rho_\omega),\,\nu\bigr)}
\end{multline}
  which takes symmetric weights $\nu\in L_\omega$ as arguments. Note that from
  the property $(\omega x,y)=(x,\omega^{-1} y)$ and the definition of
  $\proj$
  it follows $(\proj j,\nu)=(j,\nu)$ for $\nu\in L_\omega$.

\begin{lemma}
  \label{aplm:faux}
  The function $f$ is invariant under the action of
  $W_\omega\ltimes(k+g^\vee)(L^\vee)_\omega$ and vanishes for elements
  on the boundary of the Weyl chambers, in particular on
  $S_{k+g^\vee}^+\backslash S_{k+g^\vee}^{++}$.
\begin{proof}
  The property $f(w\nu)=f(\nu)$ for $w\in W_\omega$ is proved by using
  $(wx,y)=(x,w^{-1}y)$, invariance of the weight system $M_i$
  under Weyl transformations and redefinition
  of $j,w_1,w_2$. Due to $\epsilon_\omega(w)^2=1$ possible signs cancel.
  As $\proj j+w_1(\alpha+\rho_\omega)-w_2(\beta+\rho_\omega)\in
  \bigl((L^\vee)_\omega\bigr)^\vee$, the property
  $f\bigl(\nu+(k+g^\vee)\beta\bigr)=f(\nu)$ for
  $\beta\in(L^\vee)_\omega$ is obvious.
  To prove the second statement let us define the
  auxiliary function $g(\nu)=S_{\alpha,\nu-\rho}^\omega$ which enters
  each summand of the function $f(\nu)$ as a factor.
  Similar as for $f(\nu)$ one shows that
  $g\bigl(w\nu+(k+g^\vee)\beta\bigr)=\epsilon_\omega(w)g(\nu)$
  for all $\beta\in(L^\vee)_\omega$ and $w\in W_\omega$. Let $\nu$ be an element of the boundary of the
  fundamental Weyl chamber, i.e. $\nu\in S_{k+g^\vee}^+\backslash S_{k+g^\vee}^{++}$.
  Then it is either invariant under an elementary reflection or a combined
  action of a translation and an elementary reflection $w\in W_\omega$. The
  equation $g(\nu)=\epsilon_\omega(w)g(\nu)$ now implies
  that $g(\nu)=0$ and thus $f(\nu)=0$ for
  $\nu\in S_{k+g^\vee}^+\backslash S_{k+g^\vee}^{++}$.
\end{proof}
\end{lemma}

\begin{corollary}
  \label{apcl:faux}
  Eq.~\eqref{apeq:AnnulusAux} can be rewritten as
\begin{eqnarray}
  \label{apeq:faux2}
  N_{i\alpha}^\beta
  =\frac{1}{|W_\omega|}\sum_{w\in W_\omega}\sum_{\nu\in S_{k+g^\vee}^+}f(w\nu)
  =\frac{1}{|W_\omega|}\sum_{\nu\in L_\omega/(k+g^\vee)(L^\vee)_\omega}f(\nu)\:.
\end{eqnarray}
\end{corollary}
   
\begin{lemma}
  \label{aplm:aux}
  Let $\Gamma$ be a lattice and $\Gamma_s\subset\Gamma$ be a sublattice
  of the same rank as $\Gamma$. Let
  $\Gamma^\vee$ and $\Gamma_s^\vee$ be the dual lattices to
  $\Gamma,\Gamma_s$ with respect to an inner product
  $(\cdot,\cdot)$. For any $h\in\Natural$ and
  $x\in\Gamma_s^\vee$ we have
\begin{eqnarray*}
  \sum_{y\in \Gamma/h\Gamma_s}e^{2\pi i(x,y)/h}
  =\bigl|\Gamma/h\Gamma_s\bigr|\cdot\delta_{x\in h\Gamma^\vee}\:.
\end{eqnarray*}
\begin{proof}
We will use the fact that the characters $\chi$ of irreducible representations
of a finite group $G$ are orthogonal in the sense that 
$\sum_{g\in G}\chi(g)\overline{\chi^\prime(g)}
 =|G|\cdot\delta_{\chi,\chi^\prime}$.
The quotient $\Gamma/h\Gamma_s$ is a finite abelian group. For 
$x\in\Gamma_s^\vee$ the function $\chi_x:\Gamma/h\Gamma_s\to\Complex,\:
\chi_x(y)=e^{2\pi i(x,y)/h}$ is the character of an irreducible representation
of $\Gamma/h\Gamma_s$ and the character $\chi_0$ of the trivial representation
is identical to one. The orthogonality relation reads, for $x\in\Gamma_s^\vee$,
\begin{eqnarray*}
  \sum_{y\in\Gamma/h\Gamma_s}e^{2\pi i(x,y)/h}
  =\sum_{y\in\Gamma/h\Gamma_s}\chi_x(y)\overline{\chi_0(y)}
  =\bigl|\Gamma/h\Gamma_s\bigr|\cdot\delta_{\chi_x,\chi_0}\:.
\end{eqnarray*}
  But $\chi_x\equiv\chi_0$ is equivalent to $x\in h\Gamma^\vee$.
\end{proof}
\end{lemma}

\begin{proof}[Proof of Theorem~\ref{apth:main}]
  We insert expression~\eqref{apeq:fAux} for $f(\nu)$ into~\eqref{apeq:faux2}
  and apply Lemma~\ref{aplm:aux} with $\Gamma=L_\omega$,
  $\Gamma_s=(L^\vee)_\omega$
  and $h=k+g^\vee$. This results in
\begin{eqnarray*}
  N_{i\alpha}^\beta
  =\frac{1}{|W_\omega|}\sum_{j\in M_i}\sum_{w_1,w_2\in W_\omega}\!\epsilon_\omega(w_1w_2)\:\delta_{\proj j+w_1(\alpha+\rho_\omega)-w_2(\beta+\rho_\omega)\in(k+g^\vee)(L_\omega)^\vee}\:.
\end{eqnarray*}
  Using the invariance of all quantities under $W_\omega$ we are lead to the
  final result
\begin{eqnarray}
  \label{apeq:last}
  N_{i\alpha}^\beta
  =\sum_{j\in M_i}\sum_{w\in W_\omega}\epsilon_\omega(w)
  \;\delta_{w(\proj j+\alpha+\rho_\omega)-(\beta+\rho_\omega)
  \,\in\, (k+g^\vee)(L_\omega)^\vee}\:.
\end{eqnarray}
  The interpretation of the last formula then amounts to
  the algorithm of the theorem. Step 5 follows since $\beta+\rho_\omega$ is always
  in $B_{k+g^\vee}^{++}$.
\end{proof}
  
Note that for inner automorphisms the sets $P_k^+$, $S_k^+$ and $B_k^+$
all coincide and we recover the Kac-Walton formula for ordinary fusion rules.
Formula~\eqref{apeq:last} directly shows that the twisted fusion rules are
integer numbers but does not show that they are non-negative. (However,
non-negativity follows from the general theory of symmetry breaking boundary
conditions \cite{Fuchs:1999zi}.) We have also implemented the algorithm on a 
computer and have verified that no negative integers appear for the first few 
levels in the cases listed in Table~\ref{aptb:THR}.

\medskip
  
{\bf Acknowledgements} --- We would like to thank S.\ Fredenhagen and
J.\ Fuchs for a careful reading of the manuscript. C.S.\ and I.R.\ 
are grateful to the 
Albert-Einstein-Institute in Golm for hospitality. I.R.\ is
supported by the European Commission contract $\#$HPMF-CT-2000-00747.
The work of T.Q.\ is supported by the Studienstiftung des deutschen Volkes.

\providecommand{\href}[2]{#2}\begingroup\raggedright\endgroup

\end{document}